\documentclass[10pt]{amsart}
\usepackage{amsmath,amssymb,enumerate}

\newtheoremstyle{plainsl}%
	{\topsep}
	{\topsep}
	{\slshape} 
	{}
	{\normalfont\bfseries}
	{.}
	{ }
	{}
\swapnumbers

\newtheorem{theorem}{Theorem}[section]
\newtheorem{lemma}[theorem]{Lemma}

\newcommand\cref[1]{Corollary~\ref{cor:#1}}

\renewcommand\proof{\noindent\textsl{Proof. }}
\newcommand\sqr[2]{{\vbox{\hrule height.#2pt
    \hbox{\vrule width.#2pt height#1pt \kern#1pt
        \vrule width.#2pt}\hrule height.#2pt}}}
\renewcommand\qed{%
	\ifmmode\eqno\sqr53
	\else\nolinebreak\ \hfill\sqr53\medbreak\fi}

\title{Colouring the Sphere}
\author{C.~D.~G\lowercase{odsil} \\
J.~Z\lowercase{aks} $^1$    \\
                  {D}\lowercase{epartment of} C\lowercase{ombinatorics and} 
				O\lowercase{ptimization}\\
                  U\lowercase{niversity of} W\lowercase{aterloo},
                  W\lowercase{aterloo}, O\lowercase{ntario} \\
                  C\lowercase{anada} N2L 3G1}

\thanks{$^1$Usual address: Department of Mathematics, University of Haifa, Haifa, Israel. 
The work of both authors was supported in part by NSERC.}

\begin{document}

\maketitle

\begin{abstract}
Let $G$ be the graph with the points of the unit sphere in $\mathbb{R}^3$ as its vertices, by
defining two unit vectors to be adjacent if they are orthogonal as vectors. We show that the
chromatic number of this graph is four, thereby answering a question of Peter Frankl. We also
prove that the subgraph of G induced by the unit vectors with rational coordinates is
3-colourable.
\end{abstract}

\section{Introduction}
\label{sec:intro}

Let $G$ be the graph with the points of the unit sphere in three dimensions as its vertices,
and two vertices are adjacent if and only if they are orthogonal as unit vectors.

\begin{lemma}
The chromatic number of $G$ is four.
\end{lemma}

\proof 
It is clear that $G$ has chromatic number at least three. Suppose that we have a
3-colouring of $G$ using the colors red, white and black. If $x$ is a vertex in $G$ we may
view it as the north pole. If $x$ is coloured red then the points on the equator with respect
to $x$ must be coloured black and white and, furthermore, points of both colours must occur
(since any great circle on the sphere contains orthogonal pairs of vectors). It follows in
turn that the South pole with respect to $x$ must also be coloured red. Thus we have shown
that in a 3-colouring of $G$, any antipodal pair of points must have the same colour. This
means that any 3-colouring of $G$ determines a 3-colouring of the real projective plane, with
the property that the points on any line are coloured with at most two colours.

But it has been shown \cite{proj_color} that the 3-colourings of the projective plane over a
field $K$ are in one-to-one correspondence with the non-Archimedean valuations of $K$. A
non-Archimedean valuation on a field $K$ is a real-valued function $\nu$ on $K$ such that:
\begin{enumerate} 
	\item $\nu(x) \ge 0$ for all $x$ in $K$, and $\nu(x) = 0$ if and only if $x=0$. 
	\item $\nu(xy) = \nu(x)\nu(y)$ for all $x$ and $y$ in $K$. 
	\item $\nu(x+y) \le \max\{\nu(x),\nu(y)\}.$ 
\end{enumerate} 
The colouring corresponding to a non-Archimedean
valuation $\nu$ is described as follows. A point with homogeneous coordinates $(x,y,z)$ in the
projective plane over $K$ is coloured 
\begin{align*} 
	\text{red if:}& \ \ \ \ \nu(x) > \nu(y), \ \ \nu(x) > \nu(z), \\ 
	\text{white if:}& \ \ \ \ \nu(x) \le \nu(y), \ \ \nu(y) > \nu(z), \\
	\text{red if:}& \ \ \ \ \nu(x) \le \nu(z), \ \ \nu(y) \le \nu(z). 
\end{align*}

If $L$ is a subfield of $K$ then $\nu$ is also a valuation on $L$. Conversely, if a valuation
on $L$ is given then it can be extended to a valuation on $K$. By way of (a pertinent)
example, we describe a non-Archimedean valuation on the rationals $\mathbb{Q}$. Choose a prime
$p$ (any prime). Every rational number $b$ can now be written in the form $ap^n$, where the
numerator and denominator of $a$ are not divisible by $p$. Define $\nu(b)$ to be $2^{-n}$;
this is the required valuation. This is known as the $p$-adic valuation on $\mathbb{Q}$. Every
non-Archimedean valuation on $\mathbb{Q}$ is $p$-adic for some prime $p$. (For more
information on valuations see, e.g., \cite{local_fields,algebra}.)

Since the real projective plane can be 3-coloured, it follows that the sphere in
$\mathbb{R}^3$ can be 3-coloured in such a way that antipodal pairs of points have the same
colour, and all great circles are 2-coloured. But we need a colouring with the property that
orthogonal vectors have different colours.

We will call a vector $x$ in the unit sphere a $\mathbb{Q}$-point if the line it spans
contains a point with coordinates in $\mathbb{Q}$. Any 3-colouring of the real projective
plane determines a 3-colouring of the rational prjective plane , and hence a 3-colouring of
the $\mathbb{Q}$-points on the unit sphere. As noted above, such a colouring comes from one of
the $p$-adic valuations, and the question now is whether any of these give rise to a
3-colouring of the $\mathbb{Q}$-points on the sphere satisfying the orthogonality condition.
But if $p$ is odd then the points $(1,1,0)$ and $(-1,1,0)$ in the projective plane will be
coloured white. These two points determine orthogonal points on the sphere. If $p$ is equal to
2 then we see that the points $(-7,0,1)$ and $(1,0,7)$ are both coloured black, and again
determine orthogonal points on the sphere. Hence no 3-colouring of the sphere works. (We could
also have observed that the points $(-1,3,1)$ and $(-5,2,1)$ will both be coloured black, for
any prime $p$, but the argument used actually shows that the subgraph of $G$ induced by the
points in $\mathbb{Q}\left(\sqrt{2}\right)$ is not 3-colourable.)

To complete the proof, we give a four colouring of $G$. Begin by colouring the two points on
the $x$-axis red, the points on the $y$-axis white and those on the $z$-axis black. There are
exactly three great circles which pass through two of these three pairs of points. We colour
these using the two colours used on the four points that they pass through. These circles
divide the sphere into eight open octants. Colour the four octants in the half space $z > 0$
red, white, black and blue in some order. Give the remaining four octants the colour of their
antipodal octant. The reader is invited to check that this works. \qed

In contrast to the above result, we have the following:

\begin{lemma}
The chromatic number of the subgraph of $G$ induced by the unit vectors with all coordinates
rational is three.
\end{lemma}

\proof
To show this, we use the 3-colouring of the rational projective plane corresponding to the
2-adic valuation on $\mathbb{Q}$. This plane can be coordinatized by the integers. If we do
this then $(x,y,z)$ is coloured:
\begin{align*}
\text{red if:}& \ \ \ x \text{ is odd, } y \text{ and } z \text{ are even,} \\
\text{white if:}& \ \ \ y \text{ is odd and } z \text{ is even,} \\
\text{black if:}& \ \ \ z \text{ is odd.}
\end{align*}
Note that we assume that $x$, $y$ and $z$ have no common factor, and so one of these three
integers must be odd.

The point $(x,y,z)$ in the rational projective plane determines a rational point on the sphere
if and only if $x^2 + y^2 + z^2$ is a square. A necessary condition for this to occur is that
precisely one of $x$, $y$ and $z$ be odd. (This follows from the observation that any square
is congruent to 0 or 1, modulo 4, and from the fact that at least one of $x$, $y$ and $z$ is
odd.) It follows that if $(x,y,z)$ is coloured white then $x$ is also even, and if it is
coloured black then $x$ and $y$ must be even. Thus, if $(x,y,z)$ and $(x',y',z')$ are two
points in the rational projective plane with the same colour and if these points correspond to
rational points on the sphere, the inner product $xx' + yy' + zz'$ must be odd. Hence it is
not zero, and so the subgraph of $G$ induced by the rational points of the sphere is
3-chromatic.\qed

Each colour class in the above 3-colouring is dense in the sphere. To prove this, let $\alpha$
be chosen so that $\sin \alpha = \frac{3}{5}$ and $\cos \alpha = \frac{4}{5}$. Then $\alpha$
is not a rational multiple of $\pi$, and therefore $\sin n\alpha$ and $\cos n\alpha$ are
non-zero for all integers $n$. Let $F$ be the matrix which represents rotation about the
$z$-axis through an angle $\alpha$. It follows that that image $\mathcal{I}$, under the powers
of $F$, of the point $(1,0,0)$ is a dense subset of the equator. Suppose the point \[u :=
\left(\frac{a}{d},\frac{b}{d},0\right)\] is on the unit sphere and that $a$ and $d$ are odd
and $b$ is even. Then $u$ has the same colour as $(1,0,0)$, and so does $Fu$. This proves that
$\mathcal{I}$ is monochromatic. Now let $G$ represent rotation through an angle of $\alpha$
about the $y$-axis. Then the same arguments show that the image of $\mathcal{I}$ under the
powers of $G$ is a dense monochromatic set of the rational points on the unit sphere. Since
the map which sends $(x,y,z)$ to $(y,z,x)$ permutes the three colour classes cyclicly, it
follows that all three colour classes are dense. (This also shows that the rational points on
the unit sphere in $\mathbb{R}^3$ are dense. The easiest way to verify this is to note that,
under stereographic projection, the rational points on the sphere, excepting the North pole,
are in one-to-one correspondence with the rational points in the affine plane. The key to this
is that stereographic projection is a birational mapping. See e.g., Chapter XI, in particular
Exercise 23, of \cite{alg_proj_geom}.)

Despite the fact that the colour classes of this colouring are dense, there are many
monochromatic circles. To see this, choose two unit vectors $u$ and $v$ with the same colour.
Then, if $u \ne -v$, the unit vectors $x$ such that $x \cdot u = v \cdot u$ form a circle. It
is easy to verify that if \[u = \frac{1}{d}(a,b,c), \ \ \ \ \ v = \frac{1}{d'}(a',b',c')\]
then $(aa' + bb' + cc')$ is odd, and hence that the given circle is monochromatic. (If $v$ has
a different colour from $u$ then this circle will contain only rational points with a
different colour from $u$, since in this case $(aa' + bb' + cc')$ is even.)

Since $\chi(G) = 4$ it follows from a well known result of Erd\H{o}s and Szekeres that $G$
contains finite 4-chromatic subgraphs, but we have not yet found one. Note that if we have a
set of $n$ points in the sphere inducing a finite 4-chromatic subgraph $H$ of $G$, the points
in this set cannot all have rational coordinates. (Because the rational points can be
3-coloured and $H$ cannot be.)

\section*{Note}

This article first appeared as research report CORR 88-12 from the Department of Combinatorics
and Optimization at the University of Waterloo (dated March 1988). David Roberson transcribed
this into latex. Apart from a few small corrections (and this note), this version is as 
close to the original as a reasonable amount of effort could make it.

The report was never published because the main result (Lemma 1.1) was already known. However
the proof of Lemma 1.1 is new and Lemma 1.2 is also new. Since this report was widely cited
but not widely available, David felt that it would be worth posting a version on the arxiv. 
The authors thank him for
his effort.

\renewcommand{\bibname}{References}


\end{document}